\documentclass[a4paper,12pt]{article}
\usepackage{amssymb}
 \usepackage{amsthm}
\usepackage{latexsym}
\usepackage{amsmath}
\usepackage{amssymb}

\newcommand{\udots}{\mathinner{\mskip1mu\raise1pt\vbox{\kern7pt\hbox{.}}
\mskip2mu\raise4pt\hbox{.}\mskip2mu\raise7pt\hbox{.}\mskip1mu}}
\usepackage{amsmath}

\newtheorem{theorem}{Theorem}[section]
\newtheorem{corollary}{Corollary}[section]
\newtheorem{lemma}{Lemma}[section]
\newtheorem{remark}{Remark}[section]

\begin{document}
\title{\bf \Large {\bf    Hessian and increasing-Hessian orderings of multivariate skew-elliptical
random vectors  }}
{{\author{\normalsize{Chuancun Yin}\\{\normalsize\it    (School of Statistics and Data Science, Qufu Normal University}\\
\noindent{\normalsize\it Shandong 273165, China)}\\
}
\maketitle
\vskip0.01cm
\noindent{\large {\bf Abstract}}  In this work, we establish some stochastic comparison results for   multivariate skew-elliptical
random vectors. These multivariate stochastic comparisons involve Hessian and increasing-Hessian orderings as well as many of their
special cases. Necessary and/or sufficient conditions of the orderings are
provided simply based on a comparison of the underlying model parameters.
   }


\baselineskip =20pt

\numberwithin{equation}{section}
\section{Introduction}\label{intro}
The theory of stochastic orders provides an useful tool for    comparing random variables and vectors. Nowadays, stochastic orders are applied in numerous fields like  actuarial science, risk management, economics, decision theory, reliability, quality control, medicine and other related fields.  For more details, see  the monographs of Denuit et al. (2005), M\"uller and Stoyan (2002) and Shaked and Shanthikumar (2007). Among the numerous  stochastic orders, Hessian orders and increasing-Hessian orders have attracted more and more attention recently. Hessian orderings include many well-known multivariate stochastic orderings such as convex, componentwise convex, directionally-convex, supermodular, copositive, and completely-positive orderings. Increasing-Hessian orderings consist of some important orderings such as multivariate usual, increasing-convex, increasing-supermodular and increasing-completely-positive orderings; see Arlotto and Scarsini (2006),  Amiri and
Balakrishnan (2022) and  Amiri,  Balakrishnan and Eftekharian (2022).

 Hessian orders and increasing-Hessian orders of multivariate normal and elliptical distributions have been discussed rather extensively
in the literature. For instance,  results on the stochastic ordering of the multivariate normal distributions can
be found in      Houdr\'e et al. (1998),  M\"uller (2001) and Denuit and M\"uller  (2002).
 Stochastic ordering of the multivariate elliptical distributions was introduced by
Landsman and Tsanakas (2006), Davidov and Peddada (2013) and Pan, Qiu, and
Hu (2016). Yin (2021) extended the result   to general multivariate elliptical distributions.   Amiri, Izadkhah and Jamalizadeh (2020),
 Amiri and Balakrishnan (2022) considered   the    linear orderings,     Hessian and increasing-Hessian orderings of scale-shape
mixtures of multivariate skew-normal distributions. Amiri and   Balakrishnan (2022) established  some stochastic comparison results for the class of scale-shape
mixtures of multivariate skew-normal  distributions.   A unified study to a class of generalized location-scale mixture
of multivariate elliptical distributions can be found in Pu, Zhang and Yin (2023).
Related results  for  matrix variate skew-normal distributions can be found in Pu, Balakrishnan and Yin (2022).
All above discussed are restrictive on  the elliptical distributions, generalized location-scale mixture
of multivariate elliptical distributions, skew-normal distributions and the  scale-shape
mixtures of multivariate skew-normal distributions.  However, stochastic orderings of the multivariate skew-elliptical distributions have   not been well studied  yet. It remains an open problem whether necessary and/or sufficient conditions exist for stochastic orderings of multivariate skew-elliptical distributions.
This paper intends to consider these issues.

The rest of this paper is organized as follows. In   Section 2, we review some results on   multivariate skew-elliptical
distributions and the notion of stochastic orderings. In addition,  we establish  an identity for   $\mathbb{E}[f({\bf Y})]- \mathbb{E}[f({\bf X})]$, where the function $f$ fulfills some weak regularity conditions. Section 3 presents some general results for Hessian,  increasing-Hessian orderings and    linear orderings and
their special cases. Finally, some concluding remarks are made in Section 4.
\section{Preliminaries}

\subsection{ Skew-elliptical distributions}

 In recent years there has been considerable interest in the subject of   skew-elliptical distributions which allow
the handling of asymmetry and tail weight behavior simultaneously. An initial formulation of   skew-elliptical distributions  has been mentioned  in   Azzalini and  Capitanio (1999), Branco and Dey (2001)  considered the construction leading to the multivariate  skew-normal distribution  based on the conditioning mechanism. Azzalini and Capitanio (2003) examines  the skew-elliptical distributions,  they show that, similarly to the  skew-normal  family, the skew-elliptical distributions can be obtained by an additive construction approach. The developments related to multivariate skew-normal and  skew-ellptical distributions are well-summarized in Genton (2004), Lee and McLachlan (2022), Azzalini (2022).
An $n$-dimensional random variable ${\bf X}$ is said to have a multivariate skew-elliptical  distribution, if its density has the following form
\begin{equation}
f({\bf x})=2| \mathbf{ \Omega}|^{-\frac{1}{2}}g(({\bf x}- {\boldsymbol \mu})' \mathbf{ \Omega}^{-1}({\bf x-{\boldsymbol \mu} }))
 F_{g}( {\boldsymbol \alpha}'({\bf x-{\boldsymbol \mu} })), \;{\bf x}\in \Bbb{R}^n,
\end{equation}
where
$$| \mathbf{ \Omega}|^{-\frac{1}{2}}g(({\bf x}- {\boldsymbol \mu})' \mathbf{ \Omega}^{-1}({\bf x-{\boldsymbol \mu} }))$$
is a probability density function  of $n$-variate elliptical distribution ${\bf Y} \sim E_n ({\boldsymbol \mu}, \mathbf{ \Omega}, g)$,   $\mathbf{ \Omega}=(\omega_{ij}) $ is  a full rank $n\times n$   scale matrix associated to the vector ${\bf Y}$, $F_{g}(\cdot)$ is a univariate cdf of standard elliptical distribution with density generator $g$. We denote it by  ${\bf X}\sim SE_{n}\left(\boldsymbol{\mu}, \mathbf{ \Omega},\boldsymbol{\alpha}, g\right)$. For more details, see  Branco and Dey (2001).
The skew-elliptical  distributions   include the more familiar
skew-normal, skew-$t$ and skew-Chaucy distributions.

Denoted by   ${\bf \omega}={\rm diag}\{\sqrt{\omega_{11}},\cdots,\sqrt{\omega_{nn}}\},$ then, ${\bf X}$ admits the stochastic representation ${\bf X}=\boldsymbol{\mu}+{\bf \omega}{\bf Y}$, where
 \begin{eqnarray}
 {\bf Y}=\boldsymbol{\delta} |U_0|+  \mathbf{\Delta}{\bf U_1}.
  \end{eqnarray}
Here,   $\boldsymbol{\delta}=(\delta_1,\cdots,\delta_n)'$ with $\delta_1,\cdots,\delta_n\in (-1,1)$,
and   $$\left(U_{0},\mathbf{U}_1'\right)'\sim EC_{n+1}\left( \left(0,\boldsymbol{0}'\right)',\mathbf{\Psi}^*,g\right),$$
 with $\boldsymbol{0}=(0_{1},0_{2},...,0_{n})'$,  and
 \begin{equation*}
\mathbf{\Psi}^* =
\left(
\begin{array}{ccc}
  1 &  \boldsymbol{0}'  \\
 \boldsymbol{0}  & \mathbf{\Psi}- \boldsymbol{\lambda}' \boldsymbol{\lambda}
\end{array}
\right),
 \end{equation*}
is the   correlation matrix. Here,  $\boldsymbol{\lambda} =\left(\lambda_{1},\lambda_{2},...,\lambda_{n}\right)'$
with
${\lambda}_i=\delta_i(1-\delta_i^2)^{-\frac12}, \;i=1,\cdots,n$.

 The characteristic function of multivariate skew-elliptical vector ${\bf X}$ is  as follows (see Yin and Balakrishnan  (2023)).
\begin{lemma}
Assume   ${\bf X}\sim SE_{n}\left(\boldsymbol{\mu}, \mathbf{ \Omega},\boldsymbol{\alpha}_w, g^{(n+1)}\right)$ has pdf (2.1),
then its characteristic function  is given by
\begin{eqnarray}
E(e^{i{\bf t}'\mathbf{X}})=2e^{i{\bf t}'\boldsymbol{\mu}} \int_{0}^{\infty} \left(e^{i({\bf t}'\boldsymbol{\delta}_w)u_0}\cdot \phi_{L|u_0}\left({\bf t}'  \mathbf{\Omega} {\bf t}\right)\right)P(U_0\in du_0),
\end{eqnarray}
where  $\phi_{L|u_0}(\cdot)$ is a   characteristic generator and $\boldsymbol{\delta}_\omega={\bf  \omega}\boldsymbol{\delta}$.
\end{lemma}

\subsection{ Stochastic orderings}

If $f: \Bbb{R}^n \rightarrow \Bbb{R}$ is twice continuously differentiable, we write as usual
$$\nabla_f({\bf x})=\left(\frac{\partial }{\partial x_1}f({\bf x)},\cdots, \frac{\partial }{\partial x_n}f({\bf x)}\right)',\;
 {\bf H}_f({\bf x}))=\left(\frac{\partial^2}{\partial x_i \partial x_j}f({\bf x)} \right)_{n\times n}$$
for the gradient and the Hessian matrix of $f$.

A function is supermodular if and only if its Hessian matrix has nonnegative off-diagonal elements, i.e.
  $f$ is supermodular if and only if $\frac{\partial^2}{\partial x_i \partial x_j}f({\bf x)}\ge 0$ for every $i\neq j$ and ${\bf x}\in {\Bbb{R}}^n$. $f$ is said to be increasing if and only if $\nabla_f({\bf x})\ge 0$ for all ${\bf x}\in \Bbb{R}^n$; $f$ is said to
be convex if, and only if,  ${\bf H}_f({\bf x})$ is positive semi-definite for all ${\bf x}\in \Bbb{R}^n$;  $f$ is directionally convex if and only if $\frac{\partial^2}{\partial x_i \partial x_j}f({\bf x)}\ge 0$ for $1\le i, j\le n$ and ${\bf x}\in {\Bbb{R}}^n$.

An $n\times n$ matrix ${\bf A}$ is called copositive if the quadratic form ${\bf x'Ax} \ge 0$ for all ${\bf x}\ge 0$, and ${\bf A}$ is called completely positive if there exists a nonnegative $m\times n$  matrix ${\bf B}$ such that ${\bf A}={\bf B'B}$.

Denote by ${\cal C}_{cop}$ the cone of copositive matrices and by ${\cal C}_{cp}$ the cone of completely positive matrices. Let  ${\cal C}^*_{cop}$ and  ${\cal C}^*_{cp}$ be the dual of ${\cal C}_{cop}$ and  ${\cal C}_{cp}$, respectively. It is well known that (see Theorem 16.2.1 in
Hall (1986)) both cones ${\cal C}_{cp}$ and  ${\cal C}_{cop}$ are  closed and convex. Moreover,
 ${\cal C}^*_{cop}={\cal C}_{cp}$  and ${\cal C}^*_{cp}={\cal C}_{cop}$.

An integral stochastic order is a stochastic order which can be characterized by means of the comparison of their expectations.  Let $\cal{F}$ be some class of measurable functions $f: \Bbb{R}^n \rightarrow \Bbb{R}$. A stochastic order ${\bf X}\le _{\cal{F}} {\bf Y}$ is said to be integral when there exists a set $F$ of real measurable
mappings such that for two random vectors ${\bf X}$  and ${\bf Y}$,  if
 \begin{equation}
 E[f({\bf X})]\le E[f({\bf Y})]\;\; \forall f\in \cal{F}
 \end{equation}
  whenever the expectations are well defined.
A unified treatment of this type of orderings can be found in M\"uller (1997). Denuit and M\"uller (2002) proved that these orders can be defined by considering
only infinitely differentiable functions.

{\bf Definition 2.1.} A subset ${\cal C}$ of a vector space $V$ is a cone if, $x\in  {\cal C}$ and for every $\lambda\ge 0$, $\lambda x\in  {\cal C}$. The cone   ${\cal C}$ is convex if and only if     $\alpha x+\beta y\in {\cal C}$ and for every $x, y\in {\cal C}$ and $\alpha,\beta\ge 0$.  In addition, if  ${\cal C}$ is closed under an inner product
$\langle,\rangle$, then ${\cal C}^*=\{y\in V: \langle x,y\rangle\ge0, \forall x\in {\cal C}\}$ is said to be the dual of  ${\cal C}$. Moreover,
  ${\cal C}$ is said to be self-dual whenever  ${\cal C}^*={\cal C}$.

Let ${\Bbb S}$ be the space of symmetric $n\times n$  matrices and ${\Bbb {H}}$ be a closed convex cone in  ${\Bbb S}$  with inner product
 $\langle A, B\rangle =tr(AB)$ for $A, B\in  {\Bbb S}$. Then, we define the class of functions
 $${\cal{F}_{\Bbb{H}}}=\{f: {\Bbb {R}}^n \rightarrow {\Bbb{R}}: {\bf H}_f({\bf x})\in {\Bbb {H}},\;\forall x\in {\Bbb {R}}^n \}, $$
 and
  $${\Bbb{I}}=\{f: {\Bbb {R}}^n \rightarrow {\Bbb{R}}: \nabla_f({\bf x})\ge 0,\;\forall x\in {\Bbb {R}}^n \}. $$

Let ${\cal C}_{psd}$ and ${\cal C}_+$ be the cones of  positive semi-definite and non-negative matrices in  ${\Bbb S}$, respectively. Then,
these cones are closed and convex with corresponding duals   ${\cal C}^*_{psd}={\cal C}_{psd}$ and ${\cal C}^*_+={\cal C}_+$; see Theorem
16.2.1 of Hall (1986).
  Let ${\cal C}_{+off}$ and ${\cal C}_{+diag}$ be the cones of   non-negative off-diagonal elements and non-negative
main diagonal elements matrices, respectively.  Then,
these cones are closed and convex with corresponding duals (see Arlotto and Scarsini (2009))  ${\cal C}^*_{+off}$ and ${\cal C}^*_{+diag}$, where
    $${\cal C}^*_{+off}=\{C\in  {\Bbb S}: c_{ii}=0, c_{ij}\ge 0,\; for\; i\neq j\in\{1,\cdots,n\} \},$$
and
 $${\cal C}^*_{+diag}=\{C\in  {\Bbb S}: c_{ii}\ge 0, c_{ij}= 0,\; for\; i\neq j\in\{1,\cdots,n\} \}.$$
 Then, $f$ is convex if and only if (iff )  ${\bf H}_f({\bf x})\in {\cal C}_{psd}$, directionally-convex iff ${\bf H}_f({\bf x})\in {\cal C}_+$, componentwise-convex
iff ${\bf H}_f({\bf x})\in  {\cal C}_{+diag}$, and supermodular iff ${\bf H}_f({\bf x})\in  {\cal C}_{+off}$.

We now introduce the Hessian and increasing-Hessian orderings, as special integral orderings, as follows.

{\bf Definition 2.2} (Arlotto and Scarsini (2009))  If ${\cal{F}}= {\cal {F}_{\Bbb{H}}}$ in (2.1), then the order is called Hessian order. If
 ${\cal{F}}= {\cal {F}_{\Bbb{H}}}\cap {\Bbb{I}} $ in (2.1), then the order is called  increasing-Hessian order.

 Some well-known Hessian  orderings  for random vectors are as follows.   \\
 (1)  Convex order:   ${\bf X}\le _{cx} {\bf Y}$  if  $\cal{F}$  is the class of convex functions; \\
 (2) Supermodular order: ${\bf X}\le_{sm} {\bf Y}$  if $\cal{F}$  is the class of  all twice differentiable
functions $f: \Bbb{R}^n \rightarrow \Bbb{R}$ satisfying $\frac{\partial^2 }{\partial x_i \partial x_j}f({\bf x})\ge 0$ for ${\bf x}\in \Bbb{R}^n$ and all $1\le i<j\le n$.\\
(3)  Directionally convex order:  ${\bf X}\le_{dcx} {\bf Y}$ if  $\cal{F}$  is the class of all twice differentiable
functions $f: \Bbb{R}^n \rightarrow \Bbb{R}$ satisfying    $\frac{\partial^2 }{\partial x_i \partial x_j}f({\bf x})\ge 0$ for ${\bf x}\in \Bbb{R}^n$ and all $1\le i,j\le n$.\\
(4) Componentwise convex order:  ${\bf X}\le_{ccx} {\bf Y}$  if  $\cal{F}$  is the class of all twice differentiable
functions $f: \Bbb{R}^n \rightarrow \Bbb{R}$ satisfying $\frac{\partial^2 }{\partial x_i^2}f({\bf x})\ge 0$ for ${\bf x}\in \Bbb{R}^n$ and all $1\le i\le n$.\\
(5)  Completely-positive order: $X \le _{cp} Y$ if  $\cal{F}$  is the class of all    functions $f$ such that  ${\bf H}_f({\bf x})\in {\cal C}_{cp}$.\\
(6) Copositive order: $X \le _{cop} Y$ if   $\cal{F}$  is the class of    all functions $f$ such that  ${\bf H}_f({\bf x})\in {\cal C}_{cop}$.

 Some well-known increasing-Hessian orderings  for random vectors are as follows.   \\
(a)  Usual stochastic order: ${\bf X}\le _{st} {\bf Y}$ if  $\cal{F}$ is the class of increasing functions; \\
(b) Increasing convex order:    ${\bf X}\le _{icx} {\bf Y}$  if  $\cal{F}$  is the class of  increasing convex functions; \\
(c)  Increasing supermodular order: ${\bf X}\le_{ism} {\bf Y}$  if $\cal{F}$  is the class of  all twice differentiable
    functions $f: \Bbb{R}^n \rightarrow \Bbb{R}$ satisfying  $\frac{\partial }{\partial x_i}f({\bf x})\ge 0$   and $\frac{\partial^2 }{\partial x_i \partial x_j}f({\bf x})\ge 0$ for ${\bf x}\in \Bbb{R}^n$ and all $1\le i<j\le n$.\\
(d) Increasing directionally convex:  ${\bf X}\le_{idcx} {\bf Y}$ if $\cal{F}$  is the class of  all twice differentiable
functions $f: \Bbb{R}^n \rightarrow \Bbb{R}$ satisfying $\frac{\partial }{\partial x_i}f({\bf x})\ge 0$ for ${\bf x}\in \Bbb{R}^n$ and all $1\le i\le n$, and $\frac{\partial^2 }{\partial x_i \partial x_j}f({\bf x})\ge 0$ for ${\bf x}\in \Bbb{R}^n$ and all $1\le i, j\le n$.\\
(e) Upper orthant order: ${\bf X}\le_{uo} {\bf Y}$ if  $\cal{F}$  is the class of all  infinitely  differentiable
functions $f: \Bbb{R}^n \rightarrow \Bbb{R}$ satisfying $\frac{\partial^k }{\partial x_{i_1}\cdots \partial x_{i_k}}f({\bf x})\ge 0$ for ${\bf x}\in \Bbb{R}^n$ and all $1\le {i_1}<\cdots< i_k\le n$.\\
(f) Increasing componentwise convex order: ${\bf X}\le_{iccx} {\bf Y}$   if  $\cal{F}$  is the class of all twice differentiable
functions $f: \Bbb{R}^n \rightarrow \Bbb{R}$ satisfying  $\frac{\partial }{\partial x_i}f({\bf x})\ge 0$   and  $\frac{\partial^2 }{\partial x_i^2}f({\bf x})\ge 0$ for ${\bf x}\in \Bbb{R}^n$ and all $1\le i\le n$.\\
(g) Increasing completely-positive order: $X \le _{icp} Y$ if  $\cal{F}$  is the class of all  increasing  functions $f$ such that  ${\bf H}_f({\bf x})\in {\cal C}_{cp}$.\\
(h)  Increasing copositive order: $X \le _{icop} Y$ if   $\cal{F}$  is the class of    all  increasing functions $f$ such that  ${\bf H}_f({\bf x})\in {\cal C}_{cop}$.

Finally, we present several  linear orderings in the following.\\
A function $f: {\Bbb R}^n \rightarrow {\Bbb R}$  is called linear-convex if $f({\bf x})=\psi({\bf a}'{\bf X})$ for  ${\bf a}\in  {\Bbb R}^n$ and  a scalar convex function $\psi$. A function $f: {\Bbb R}^n \rightarrow {\Bbb R}$  is called  positive-linear-convex if $f({\bf x})=\psi({\bf a}'{\bf X})$ with  ${\bf a}\in  {\Bbb R}_+^n$ and  a scalar convex function $\psi$.\\
(i) Linear-convex order: ${\bf X}\le _{lcx} {\bf Y}$ if   ${\bf a}'{\bf X}\leq_{cx}{\bf a}'{\bf Y}$, for ${\bf a}\in \Bbb{R}^n$; \\
(ii) Positive-linear-convex order:  ${\bf X}\le _{plcx} {\bf Y}$ if  ${\bf a}'{\bf X}\leq_{cx} {\bf a}'{\bf Y}$, for ${\bf a}\in \Bbb{R}_+^n$; \\
(iii) Increasing-positive-linear-convex order:  ${\bf X}\le _{iplcx} {\bf Y}$ if ${\bf a}'{\bf X}\leq_{icx} {\bf a}'{\bf Y}$, for ${\bf a}\in \Bbb{R}_+^n$; \\
(iv) Positive-linear-usual order:  ${\bf X}\le _{plst} {\bf Y}$ if   ${\bf a}'{\bf X}\leq_{st}{\bf a}'{\bf Y}$, for ${\bf a}\in \Bbb{R}_+^n$;\\
(v) Linear-usual order:  ${\bf X}\le _{lst} {\bf Y}$ if   ${\bf a}'{\bf X}\leq_{st}{\bf a}'{\bf Y}$, for ${\bf a}\in \Bbb{R}^n$.

\begin{remark} The following implications are well known (see e.g.  Shaked and Shanthikumar (2007) and  Scarsini (1998)).
$${\bf X}\le_{ilcx} {\bf Y}\Leftrightarrow {\bf X}\le_{lcx} {\bf Y},  $$
$${\bf X}\le_{cx} {\bf Y}\Leftrightarrow {\bf X}\le_{icx} {\bf Y}\; {\rm and}\; E({\bf X})=E({\bf Y}),$$
$$ {\bf X}\le_{cx} {\bf Y} \Rightarrow  {\bf X}\le_{lcx} {\bf Y} \Rightarrow  {\bf X}\le_{plcx} {\bf Y},$$
$${\bf X}\le_{st} {\bf Y}\Rightarrow  {\bf X}\le_{icx} {\bf Y} \Rightarrow  {\bf X}\le_{ilcx} {\bf Y} \Rightarrow  {\bf X}\le_{iplcx} {\bf Y}.$$
\end{remark}

\subsection{ An identity for   $\mathbb{E}[f({\bf Y})]- \mathbb{E}[f({\bf X})]$}\label{intro}

Denote by ${\bf O}_{n\times n}$ the $n \times n$ matrix  with all entries  equal to $0$  and  by ${\bf I}_n$  the $n\times n$ identity matrix. For symmetric matrices $A$ and $B$ of the same size, the notion $A\preceq B$ or $B-A  \succeq {\bf O}$ means that $B-A$ is positive semi-definite. The inequality between vectors or matrices denotes componentwise inequalities. Throughout, we assume all involved integrals
and expectations exist.

Let ${\bf X}$ and ${\bf Y}$ be $n$-dimensional random vectors from multivariate skew-elliptical family of distributions with
 \begin{equation}
 {\bf X}\sim SE_{n}\left(\boldsymbol{\mu}^x, \mathbf{ \Omega}^x,\boldsymbol{\alpha}_w^x, g^{(n+1)}\right),\;\; {\bf{Y}}\sim SE_n ({\boldsymbol \mu}^y, \mathbf{ \Omega}^y,\boldsymbol{\alpha}_w^y,  g^{(n+1)}).
 \end{equation}
The mean vector and the covariance matrix of ${\bf X}$ are, respectively,
\begin{eqnarray}
&&E({\bf X})=\boldsymbol{\mu}^x+ E(|U_0|)\boldsymbol{\delta}_w^x,\\
&&Var({\bf X})=\frac{1}{n+1}E(R_0^2)\mathbf{\Omega}^x- (E(|U_0|))^2{\boldsymbol{\delta}_w^x} {\boldsymbol{\delta}_w^x}',
\end{eqnarray}
where  $R_x$ is defined by   $\mathbf{X}\stackrel{d}{=}R_0\mathbf{M}_x+\boldsymbol{\mu}^x$; see Yin and Balakrishnan (2023).

The following lemma provides general form of  an identity for   $\mathbb{E}[f({\bf Y})]- \mathbb{E}[f({\bf X})]$, which is a key tool for studying the stochastic comparison  between ${\bf X}$ and ${\bf Y}$.

\begin{lemma} Let  ${\bf{X}}$ and ${\bf{Y}}$ be as in  (2.5).  Assume that
$f: \Bbb{R}^n \rightarrow \Bbb{R}$ is twice continuously differentiable and satisfies some polynomial
growth conditions at infinity:
\begin{equation}
f({\bf x})=O(||{\bf x}||),\; \bigtriangledown f({\bf x})=O(||{\bf x}||).
\end{equation}
 Then
\begin{eqnarray*}
  \mathbb{E}[f({\bf Y})]- \mathbb{E}[f({\bf X})]&=&\int_0^1\int_{\Bbb{R}^n}({\boldsymbol \mu}^y-{\boldsymbol \mu}^x)'\nabla f({\bf x})\phi_{\lambda}({\bf x})d{\bf x}d\lambda\\
  &&+E(|U_0|)\int_0^1\int_{\Bbb{R}^n}( {\boldsymbol \delta}_w^y-{\boldsymbol \delta}_w^x)'\nabla f({\bf x})\phi_{1,\lambda}({\bf x})d{\bf x}d\lambda\\
&&+\frac{\mathbb{E}(R^2)}{2n}E(|U_0|) \int_0^1\int_{\Bbb{R}^n}tr\{( \mathbf{\Omega}^{y}- \mathbf{\Omega}^{x}){\bf H}_f({\bf x})\}\phi_{2,\lambda}({\bf x})d{\bf x}d\lambda,
\end{eqnarray*}
where tr$(A)$ denotes the trace of the matrix $A$,   $R\ge 0$   is the random variable with $R \sim F$ in $[0, \infty)$ called the generating variate, $\phi_{\lambda}, \phi_{1,\lambda}, \phi_{2,\lambda}$ are three density functions specified below.
\end{lemma}

 {\bf  Proof.}  For $0\le \lambda\le 1$, setting
  $${\boldsymbol \mu}_{\lambda}:= \lambda{\boldsymbol \mu}^y+(1-\lambda){\boldsymbol \mu}^x,$$
  $${\boldsymbol \delta}_{\lambda}:= \lambda {\boldsymbol \delta}_w^y+(1-\lambda) {\boldsymbol \delta}_w^x,$$
  and
  $$ \mathbf{\Omega}_{\lambda}:=\lambda \mathbf{\Omega}^y+(1-\lambda) \mathbf{\Omega}^x.$$
    Define
$$ \Psi_{\lambda}({\bf t})=2 e^{i{\bf t}'\boldsymbol{\mu}_{\lambda}}\int_{0}^{\infty} \left(e^{i({\bf t}'\boldsymbol{\delta}_{\lambda})u_0}\cdot \phi_{L|u_0}
\left({\bf t}' \mathbf{\Omega}_{\lambda} {\bf t}\right)\right)P(U_0\in du_0),\; {\bf t}\in \Bbb{R}^n.$$
By using the Fourier inversion theorem
$$\phi_{\lambda}({\bf x})=\left(\frac{1}{2\pi}\right)^n\int e^{-i{\bf t}'{\bf x}}  \Psi_{\lambda}({\bf t})d{\bf t}.$$
The derivative of $\Psi_{\lambda}$ with respect to $\lambda$ is
\begin{eqnarray*}
\frac{\partial \Psi_{\lambda}({\bf t})}{\partial \lambda}&=&2i{\bf t}'({\boldsymbol \mu}^y-{\boldsymbol \mu}^x)e^{i{\bf t}'\boldsymbol{\mu}_{\lambda}}\int_{0}^{\infty} \left(e^{i({\bf t}'\boldsymbol{\delta}_{\lambda})u_0}\cdot \phi_{L|u_0}\left({\bf t}' \mathbf{\Omega}_{\lambda} {\bf t}\right)\right)P(U_0\in du_0) \\
&&+2i{\bf t}'({\boldsymbol \delta}_w^y-{\boldsymbol \delta}_w^x)e^{i{\bf t}'\boldsymbol{\mu}_{\lambda}}\int_{0}^{\infty}u_0 \left(e^{i({\bf t}'\boldsymbol{\delta}_{\lambda})u_0}\cdot \phi_{L|u_0}\left({\bf t}' \mathbf{\Omega}_{\lambda} {\bf t}\right)\right)P(U_0\in du_0)\\
&&+2{\bf t}'( \mathbf{\Omega}^y-  \mathbf{\Omega}^x){\bf t}e^{i{\bf t}'\boldsymbol{\mu}_{\lambda}}\int_{0}^{\infty} \left(e^{i({\bf t}'\boldsymbol{\delta}_{\lambda})u_0}\cdot \phi'_{L|u_0}\left({\bf t}' \mathbf{\Omega}_{\lambda} {\bf t}\right)\right)P(U_0\in du_0)\\
&=&i{\bf t}'({\boldsymbol \mu}^y-{\boldsymbol \mu}^x)\Psi_{\lambda}({\bf t})\\
&&+ iE(|U_0|)  {\bf t}' ( {\boldsymbol \delta}_w^y- {\boldsymbol \delta}_w^x)\Psi_{1,\lambda}({\bf t})\\
&&-\left(\frac{\mathbb{E}(R^2)}{2n}\right) E(|U_0|){\bf t}'( \mathbf{\Omega}^y-  \mathbf{\Omega}^x){\bf t}  \Psi_{2,\lambda}({\bf t}),
\end{eqnarray*}
where
$$\Psi_{1,\lambda}({\bf t})=\frac{e^{i{\bf t}'\boldsymbol{\mu}_{\lambda}}}{E(U_01_{\{U_0>0\}})}\int_{0}^{\infty}u_0 \left(e^{i({\bf t}'\boldsymbol{\delta}_{\lambda})u_0}\cdot \phi_{L|u_0}\left({\bf t}' \mathbf{\Omega}_{\lambda} {\bf t}\right)\right)P(U_0\in du_0), $$
and
$$\Psi_{2,\lambda}({\bf t})=\frac{e^{i{\bf t}'\boldsymbol{\mu}_{\lambda}}}{ E(U_01_{\{U_0>0\}})}\int_{0}^{\infty} \left(e^{i({\bf t}'\boldsymbol{\delta}_{\lambda})u_0}\cdot \Psi_1\left({\bf t}' \mathbf{\Omega}_{\lambda} {\bf t}\right)\right)P(U_0\in du_0), $$
which are two characteristic functions.
Here,
$$\Psi_1(u)=\frac{1}{\mathbb{E}(R^2)}\int_0^{\infty}{}_0F_1\left(\frac{n}{2}+1;-\frac{r^2 u}{4}\right)r^2\mathbb{P}(R\in dr)$$   is a characteristic generator with
$${}_0F_1(\gamma;z)=\sum_{k=0}^{\infty} \frac{\Gamma(\gamma)}{\Gamma(\gamma+k)}\frac{z^k}{k!},$$
which is the  generalized hypergeometric series of order $(0, 1)$.
Hence,
\begin{eqnarray*}
\frac{\partial \phi_{\lambda}({\bf x})}{\partial \lambda}&=&\left(\frac{1}{2\pi}\right)^n\int e^{-i{\bf t}'{\bf x}}\frac{\partial \Psi_{\lambda}({\bf t})}{\partial \lambda} d{\bf t}\\
&=&\left(\frac{1}{2\pi}\right)^n  i\int e^{-i{\bf t}'{\bf x}}\Psi_{\lambda}({\bf t}){\bf t}'({\boldsymbol \mu}^y-{\boldsymbol \mu}^x)d{\bf t}\\
&&+\left(\frac{1}{2\pi}\right)^n iE(|U_0|)\int e^{-i{\bf t}'{\bf x}}\Psi_{1,\lambda}({\bf t})   {\bf t}'({\boldsymbol \delta}_w^y- {\boldsymbol \delta}_w^x)d{\bf t}\\
&&-\left(\frac{\mathbb{E}(R^2)}{2n}\right) \left(\frac{1}{2\pi}\right)^n  E(|U_0|) \int e^{-i{\bf t}'{\bf x}}{\bf t}'( \mathbf{\Omega}^y-  \mathbf{\Omega}^x){\bf t}  \Psi_{2,\lambda}({\bf t}) d{\bf t}\\
&=&\sum_{i=1}^n(\mu_i^y-\mu_i^x)\frac{\partial \phi_{\lambda}({\bf x})}{\partial x_i}\\
&&+E(|U_0|)\sum_{i=1}^n(\delta_i^y- \delta_i^x)\frac{\partial \phi_{1,\lambda}({\bf x})}{\partial x_i} \\
&&+\frac{\mathbb{E}(R^2)}{2n} E(|U_0|)\sum_{i=1}^n\sum_{j=1}^n(\omega^y_{ij}-\omega^x_{i j})\frac{\partial^2 \phi_{2,\lambda}({\bf x})}{\partial x_i \partial x_j},
\end{eqnarray*}
where  $ \mathbf{\Omega}^x=(\omega^x_{ij})_{n\times n}$,  $ \mathbf{\Omega}^y=(\omega^y_{ij})_{n\times n}$ and
$$\phi_{k,\lambda}({\bf x})=\left(\frac{1}{2\pi}\right)^n\int e^{-i{\bf t}'{\bf x}}  \Psi_{k,\lambda}({\bf t})d{\bf t},\; k=1,2.$$
Define $g(\lambda)=\int_{\Bbb{R}^n} f({\bf x})\phi_{\lambda}({\bf x})d{\bf x}$. It follows by using integration by parts with taking into account (2.8) that
\begin{eqnarray*}
  g'(\lambda)&=&\int_{\Bbb{R}^n} f({\bf x})\frac{\partial \phi_{\lambda}({\bf x})}{\partial \lambda}d{\bf x}\\
&=&\int_{\Bbb{R}^n}({\boldsymbol \mu}^y-{\boldsymbol \mu}^x)'\nabla f({\bf x})\phi_{\lambda}({\bf x})d{\bf x}\\
&&+E(|U_0|)\int_{\Bbb{R}^n}({\boldsymbol \delta}_w^y-{\boldsymbol \delta}_w^x)'\nabla f({\bf x})\phi_{1,\lambda}({\bf x})d{\bf x}\\
&&+\frac{\mathbb{E}(R^2)}{2n} E(|U_0|) \int_{\Bbb{R}^n}tr\{( \mathbf{\Omega}^{y}- \mathbf{\Omega}^{x}){\bf H}_f({\bf x})\}\phi_{2,\lambda}({\bf x})d{\bf x}.
\end{eqnarray*} The result follows since
$\mathbb{E}[f({\bf Y})]-\mathbb{E}[f({\bf X})]=g(1)-g(0)=\int_0^1 g'(\lambda)d\lambda$.  $\hfill\square$

The following result is a direct consequence of  Lemma 2.2.

 \begin{corollary}  Let  ${\bf{X}}$ and ${\bf{Y}}$ be as in  (2.5).  Assume that
$f: \Bbb{R}^n \rightarrow \Bbb{R}$   satisfies the conditions of Lemma 2.2, and for all ${\bf x}\in \Bbb{R}^n$:
\begin{eqnarray}
&& ({\boldsymbol \mu}^y-{\boldsymbol \mu}^x)'\nabla f({\bf x})\ge 0,\\
&& ( {\boldsymbol \delta}_w^y-{\boldsymbol \delta}_w^x)'\nabla f({\bf x})\ge 0,\\
&&  tr\{( \mathbf{\Omega}^{y}- \mathbf{\Omega}^{x}){\bf H}_f({\bf x})\} \ge 0.
\end{eqnarray}
Then, we have
 $\mathbb{E} [f({\bf X})]\le\mathbb{E}[f({\bf Y})]$.
\end{corollary}

\section{Main results  }\label{intro}

In the sequel, we investigate necessary and sufficient conditions for integral stochastic orders in the   skew-elliptical   family.

\subsection{ Hessian orderings}\label{intro}

The following theorem provides general form of sufficient conditions for the Hessian   orderings of the SE family.
\begin{theorem} Let ${\Bbb {H}} \subseteq {\Bbb S}$ and let ${\cal{C}}_{\Bbb {H}}$ be a closed convex cone  generated by   ${\Bbb H}$.
Suppose the random variables  ${\bf X}$  and ${\bf Y}$  are as in (2.5).\\
If ${\boldsymbol \mu}^x= {\boldsymbol \mu}^y,   {\boldsymbol \delta}_w^x={\boldsymbol \delta}_w^y$ and
${\bf\Omega}^{y}-{\bf\Omega }^x\in {\cal{C}}_{\Bbb H}^*$, then ${\bf X}\le_{\cal {F}_{\Bbb{H}}} {\bf Y}$.
\end{theorem}
{\bf Proof}. If the conditions hold, then, by Corollary 2.1, for $f\in{\cal{F}_{\Bbb{H}}}$ we have  $\mathbb{E} [f({\bf X})]\le\mathbb{E}[f({\bf Y})]$.
 Hence, the theorem.  $\hfill\square$

In the next theorem, we present  sufficient conditions for some well known Hessian orders.
\begin{theorem} Suppose the random variables  ${\bf X}$  and ${\bf Y}$  are as in (2.5).\\
(a) If ${\boldsymbol \mu}^x= {\boldsymbol \mu}^y,   {\boldsymbol \delta}_w^x={\boldsymbol \delta}_w^y$ and ${\bf\Omega}^{y}-{\bf\Omega}^x\succeq O$, then ${\bf X}\le _{cx} {\bf Y}$;\\
(b) If ${\boldsymbol \mu}^x= {\boldsymbol \mu}^y,  {\boldsymbol \delta}_w^x={\boldsymbol \delta}_w^y,  \omega_{ii}^x= \omega_{ii}^y$ and $ \omega_{ij}^x\le  \omega _{ij}^y$ for all $1\le i\neq j\le n$, then ${\bf X}\le _{sm} {\bf Y}$;\\
(c) If ${\boldsymbol \mu}^x= {\boldsymbol \mu}^y,   {\boldsymbol \delta}_w^x={\boldsymbol \delta}_w^y$   and $\omega _{ij}^x\le \omega_{ij}^y$ for all $1\le i,j\le n$, then ${\bf X}\le _{dcx} {\bf Y}$;\\
(d) If ${\boldsymbol \mu}^x= {\boldsymbol \mu}^y,  {\boldsymbol \delta}_w^x={\boldsymbol \delta}_w^y$,   $ \omega^x_{ii}\le   \omega^y_{ii}$ for all $1\le i\le n$, and  $\omega ^x_{ij}= \omega^y_{ij}$ for all $1\le i<j\le n$, then ${\bf X}\le _{ccx} {\bf Y}$;\\
(e)  If   ${\boldsymbol \mu}^x= {\boldsymbol \mu}^y$, ${\boldsymbol \delta}_w^x={\boldsymbol \delta}_w^y$ and ${\bf\Omega}^y-{\bf  \Omega}^x$ is copositive, then ${\bf X}\le_{cp} {\bf Y}$;\\
(f)  If   ${\boldsymbol \mu}^x= {\boldsymbol \mu}^y$, ${\boldsymbol \delta}_w^x={\boldsymbol \delta}_w^y$ and ${\bf \Omega }^y-{\bf  \Omega}^x$ is  completely copositive, then ${\bf X}\le_{cop} {\bf Y}$.
\end{theorem}
{\bf Proof}\; (a).  Suppose $f: \Bbb{R}^n \rightarrow \Bbb{R}$ is a   twice differential convex function, then  ${\bf H}_f({\bf x})$ is positive semi-definite for all ${\bf x}\in \Bbb{R}^n$.  Since ${\bf\Omega}^{y}-{\bf\Omega}^x\succeq O$, we can write it as ${\bf\Omega}^{y}-{\bf\Omega}^x=\sum_{i=1}^r {\bf \beta}_i{\bf \beta}'_i$, where $r$ is the rank of ${\bf\Omega}^{y}-{\bf\Omega}^x$, ${\bf \beta}_i$ is vector with rank 1. Hence
 \begin{eqnarray*}
  tr\{( \mathbf{\Omega}^{y}- \mathbf{\Omega}^{x}){\bf H}_f({\bf x})\} =\sum_{i=1}^r {\bf \beta}'_i {\bf H}_{f({\bf x})} {\bf \beta}_i\ge 0.
\end{eqnarray*}
 By Lemma 2.1, we have that $E[f({\bf Y})]\ge E[f({\bf X})]$, which implies ${\bf X}\le _{cx} {\bf Y}$.
Assertions  (b)-(f) can be proved in the same way.  $\hfill\square$

 The following result that we present in this section gives necessary conditions for the order of two  random vectors.
 \begin{theorem} Let ${\Bbb {H}} \subseteq {\Bbb S}$ and let ${\cal{C}}_{\Bbb {H}}$ be a closed convex cone  generated by   ${\Bbb H}$.
Suppose the random variables  ${\bf X}$  and ${\bf Y}$  are as in (2.5).
If ${\bf X}\le _ {\cal {F}_{\Bbb{H}}} {\bf Y}$, then,
 $\boldsymbol{\mu}^x+ E(|U_0|)\boldsymbol{\delta}_w^x=\boldsymbol{\mu}^y+ E(|U_0|)\boldsymbol{\delta}_w^y$ and
${\bf \Omega}^{y}-{\bf  \Omega }^x\in {\cal{C}}_{\Bbb H}^*$.
\end{theorem}
{\bf Proof}. Observe that  ${\bf O}_{n\times n}\in {\cal{C}}_{\Bbb {H}}$.  For $i=1,\cdots,n$, consider the functions $f_i({\bf x})=x_i$ and $g_i({\bf x})=-x_i$. Clearly,
 $f_i, g_i \in {\cal {F}_{\Bbb{H}}}$ since ${\bf H}_{f_i}({\bf x}))={\bf H}_{f_i}({\bf x}))={\bf O}_{n\times n}$. Applying $f_i({\bf x})$ and $g_i({\bf x})$ to
${\bf X}\le _ {\cal {F}_{\Bbb{H}}} {\bf Y}$ lead to  $\boldsymbol{\mu}^x+ E(|U_0|)\boldsymbol{\delta}_w^x=\boldsymbol{\mu}^y+ E(|U_0|)\boldsymbol{\delta}_w^y$.
Let  $\boldsymbol{\mu}^x =\boldsymbol{\mu}^y = \boldsymbol{\mu}$.
Choose a matrix ${\bf A}\in {\cal{C}}_{\Bbb {H}}$  and define a function $f$ as
$$f({\bf x})=({\bf x}-\boldsymbol{\mu})'{\bf A}({\bf x}-\boldsymbol{\mu}).$$ Observe that $f\in {\cal {F}_{\Bbb{H}}}$ since
${\bf H}_{f}({\bf x}))=\frac12{\bf A}$.
Applying $f({\bf x})$ to ${\bf X}\le _ {\cal {F}_{\Bbb{H}}} {\bf Y}$ leads to
$$E({\bf X}-\boldsymbol{\mu})'{\bf A}({\bf X}-\boldsymbol{\mu})=Ef({\bf X})\le Ef({\bf Y})=E({\bf Y}-\boldsymbol{\mu})'{\bf A}({\bf Y}-\boldsymbol{\mu}).$$
This is equivalent to
$$E(tr[{\bf X}-\boldsymbol{\mu})({\bf X}-\boldsymbol{\mu})'{\bf A}])\le E(tr[{\bf Y}-\boldsymbol{\mu})({\bf Y}-\boldsymbol{\mu})'{\bf A}]),$$
which, together with (2.7), implies that
$$tr[ ({\bf \Omega}^y-{\bf \Omega }^x)'{\bf A}]\ge 0.$$
By  the arbitrariness of ${\bf A}\in {\cal{C}}_{\Bbb {H}}$,  we get ${\bf \Omega}^{y}-{\bf  \Omega }^x\in {\cal{C}}_{\Bbb H}^*$. Thus completing   the proof of Theorem 3.3.
 $\hfill\square$

In the next theorem, we present  necessary conditions for some well known Hessian orders.
\begin{theorem} Suppose the random variables  ${\bf X}$  and ${\bf Y}$  are as in (2.5).\\
(a) If ${\bf X}\le_{cx} {\bf Y}$, then,
 $\boldsymbol{\mu}^x+ E(|U_0|)\boldsymbol{\delta}_w^x=\boldsymbol{\mu}^y+ E(|U_0|)\boldsymbol{\delta}_w^y$ and ${\bf\Omega}^{y}-{\bf\Omega}^x\succeq O$;\\
(b) If ${\bf X}\le _{sm} {\bf Y}$,  then,  $\boldsymbol{\mu}^x+ E(|U_0|)\boldsymbol{\delta}_w^x=\boldsymbol{\mu}^y+ E(|U_0|)\boldsymbol{\delta}_w^y$,
$\omega_{ii}^x= \omega_{ii}^y$ and $ \omega_{ij}^x\le  \omega _{ij}^y$ for all $1\le i\neq j\le n$;\\
(c) If     ${\bf X}\le _{dcx} {\bf Y}$, then,  $\boldsymbol{\mu}^x+ E(|U_0|)\boldsymbol{\delta}_w^x=\boldsymbol{\mu}^y+ E(|U_0|)\boldsymbol{\delta}_w^y$ and    $\omega _{ij}^x\le \omega_{ij}^y$ for all $1\le i,j\le n$;\\
(d) If ${\bf X}\le _{ccx} {\bf Y}$,  then,  $\boldsymbol{\mu}^x+ E(|U_0|)\boldsymbol{\delta}_w^x=\boldsymbol{\mu}^y+ E(|U_0|)\boldsymbol{\delta}_w^y$ and    $ \omega^x_{ii}\le   \omega^y_{ii}$ for all $1\le i\le n$, and  $\omega ^x_{ij}= \omega^y_{ij}$ for all $1\le i<j\le n$;\\
(e)  If    ${\bf X}\le_{cp} {\bf Y}$, then,  $\boldsymbol{\mu}^x+ E(|U_0|)\boldsymbol{\delta}_w^x=\boldsymbol{\mu}^y+ E(|U_0|)\boldsymbol{\delta}_w^y$ and    ${\bf\Omega}^y-{\bf  \Omega}^x$ is copositive;\\
(f)  If  ${\bf X}\le_{cop} {\bf Y}$,  then,  $\boldsymbol{\mu}^x+ E(|U_0|)\boldsymbol{\delta}_w^x=\boldsymbol{\mu}^y+ E(|U_0|)\boldsymbol{\delta}_w^y$ and    ${\bf \Omega }^y-{\bf  \Omega}^x$ is  completely copositive.
\end{theorem}
{\bf Proof}.  Theorem 3.3 implies that $\boldsymbol{\mu}^x+ E(|U_0|)\boldsymbol{\delta}_w^x=\boldsymbol{\mu}^y+ E(|U_0|)\boldsymbol{\delta}_w^y$  holds in all cases (a)-(f). Now,  we only verify the conditions about  ${\bf\Omega}^{y}$ and ${\bf\Omega}^x$.

(a) If ${\bf X}\le_{cx} {\bf Y}$, then,  ${\bf\Omega}^{y}-{\bf\Omega}^x\succeq O$ since   ${\cal C}^*_{psd}={\cal C}_{psd}$.

(b) If ${\bf X}\le _{sm} {\bf Y}$,  then,
$${\bf \Omega}^{y}-{\bf  \Omega }^x\in{\cal C}^*_{+off}=\{C\in  {\Bbb S}: c_{ii}=0, c_{ij}\ge 0,\; for\; i\neq j\in\{1,\cdots,n\} \},$$
and so $\omega_{ii}^x= \omega_{ii}^y$ and $ \omega_{ij}^x\le  \omega _{ij}^y$ for all $1\le i\neq j\le n$.

(c) If     ${\bf X}\le _{dcx} {\bf Y}$, then, ${\bf \Omega}^{y}-{\bf  \Omega }^x\in {\cal C}^*_+={\cal C}_+$ and so  $\omega _{ij}^x\le \omega_{ij}^y$ for all $1\le i,j\le n$.

(d) If ${\bf X}\le _{ccx} {\bf Y}$,  then,  $$ {\bf \Omega}^{y}-{\bf  \Omega }^x\in {\cal C}^*_{+diag}=\{C\in  {\Bbb S}: c_{ii}\ge 0, c_{ij}= 0,\; for\; i\neq j\in\{1,\cdots,n\} \},$$ from which we get
 $ \omega^x_{ii}\le   \omega^y_{ii}$ for all $1\le i\le n$, and  $\omega ^x_{ij}= \omega^y_{ij}$ for all $1\le i<j\le n$.

(e) If    ${\bf X}\le_{cp} {\bf Y}$, then,  $  {\bf \Omega}^{y}-{\bf  \Omega }^x\in{\cal C}^*_{cp}={\cal C}_{cop}$, and thus  ${\bf\Omega}^y-{\bf  \Omega}^x$ is copositive.

(f)  If  ${\bf X}\le_{cop} {\bf Y}$,  then,  $  {\bf \Omega}^{y}-{\bf  \Omega }^x\in{\cal C}^*_{cop}={\cal C}_{cp}$, and thus  ${\bf\Omega}^y-{\bf  \Omega}^x$ is  completely copositive. Completing the proof of  Theorem 3.4.  $\hfill\square$

\subsection{ Increasing-Hessian orderings}\label{intro}

The following theorem provides general form of sufficient conditions for the increasing-Hessian orderings of the  skew-elliptical distributions.

\begin{theorem} Let ${\Bbb {H}} \subseteq {\Bbb S}$ and let ${\cal{C}}_{\Bbb {H}}$ be a closed convex cone  generated by   ${\Bbb H}$.
Suppose the random variables  ${\bf X}$  and ${\bf Y}$  are as in (2.5).\\
If ${\boldsymbol \mu}^x\le {\boldsymbol \mu}^y,   {\boldsymbol \delta}_w^x\le{\boldsymbol \delta}_w^y$ and
${\bf \Omega}^{y}-{\bf  \Omega }^x\in {\cal{C}}_{\Bbb H}^*$, then ${\bf X}\le _ {\cal {F}_{\Bbb{H}}\cap{\Bbb I}} {\bf Y}$.
\end{theorem}
{\bf Proof}. If the conditions hold, then, by Corollary 2.1, for $f\in{\cal{F}_{\Bbb{H}}\cap{\Bbb I}}$, we have  $\mathbb{E} [f({\bf X})]\le\mathbb{E}[f({\bf Y})]$.
Hence  the theorem.  $\hfill\square$

In the next theorem, we present some special cases of increasing-Hessian orders.
\begin{theorem} Suppose the random variables  ${\bf X}$  and ${\bf Y}$  are as in (2.5).\\
(1) If ${\boldsymbol \mu}^x\le {\boldsymbol \mu}^y, {\boldsymbol \delta}_w^x\le {\boldsymbol \delta}_w^y$ and ${\bf  \Omega}^{x} = {\bf \Omega}^y$, then ${\bf X}\le _{st} {\bf Y}$;\\
(2) If ${\boldsymbol \mu}^x\le{\boldsymbol \mu}^y,  {\boldsymbol \delta}_w^x\le {\boldsymbol \delta}_w^y$ and ${\bf  \Omega}^{y}-{\bf  \Omega}^x\succeq O$, then ${\bf X}\le _{icx} {\bf Y}$;\\
(3) If ${\boldsymbol \mu}^x\le {\boldsymbol \mu}^y,    {\boldsymbol \delta}_w^x\le {\boldsymbol \delta}_w^y,    \omega_{ii}^x=   \omega_{ii}^y$ and $   \omega_{ij}^x\le  \omega_{ij}^y$ for all $1\le i \neq j\le n$, then ${\bf X}\le _{ism} {\bf Y}$;\\
(4) If ${\boldsymbol \mu}^x\le {\boldsymbol \mu}^y, {\boldsymbol \delta}_w^x\le {\boldsymbol \delta}_w^y$  and $  \omega_{ij}^x\le   \omega_{ij}^y$ for all $1\le i,j\le n$, then ${\bf X}\le _{idcx} {\bf Y}$;\\
(5) If ${\boldsymbol \mu}^x\le {\boldsymbol \mu}^y,  {\boldsymbol \delta}_w^x\le {\boldsymbol \delta}_w^y$,   $  \omega^x_{ii}\le    \omega^y_{ii}$ for all $1\le i\le n$, and  $ \omega^x_{ij}=   \omega^y_{ij}$ for all $1\le i<j\le n$, then ${\bf X}\le _{iccx} {\bf Y}$;\\
(6)  If   ${\boldsymbol \mu}^x \le {\boldsymbol \mu}^y$, ${\boldsymbol \delta}_w^x\le {\boldsymbol \delta}_w^y$ and ${\bf  \Omega}^y-{\bf  \Omega}^x$ is copositive, then ${\bf X}\le_{icp} {\bf Y}$;\\
(7)  If   ${\boldsymbol \mu}^x \le {\boldsymbol \mu}^y$, ${\boldsymbol \delta}_w^x\le {\boldsymbol \delta}_w^y$ and ${\bf \Omega}^y-{\bf  \Omega}^x$ is  completely copositive, then ${\bf X}\le_{icop} {\bf Y}$;\\
(8) If ${\boldsymbol \mu}^x\le {\boldsymbol \mu}^y,  {\boldsymbol \delta}_w^x\le {\boldsymbol \delta}_w^y,   \omega_{ii}^x=  \omega_{ii}^y$ and $\omega_{ij}^x\le  \omega_{ij}^y$ for all $1\le i \neq j\le n$, then ${\bf X}\le _{uo} {\bf Y}$.
\end{theorem}

{\bf Proof}\; (1).  For any increasing twice differential function   $f: \Bbb{R}^n \rightarrow \Bbb{R}$. If ${\boldsymbol \mu}^x\le {\boldsymbol \mu}^y, {\boldsymbol \delta}_w^x\le {\boldsymbol \delta}_w^y$ and ${\bf  \Omega}^{x} = {\bf \Omega}^y$, then, by Lemma 2.1. we get $\mathbb{E} [f({\bf X})]\le\mathbb{E}[f({\bf Y})]$, and so  ${\bf X}\le _{st} {\bf Y}$, as required.

Assertions  (2)-(8) can be proved in the same way.  $\hfill\square$

The following result  gives necessary conditions for the order of two  random vectors.
\begin{theorem} Let ${\Bbb {H}} \subseteq {\Bbb S}$ and let ${\cal{C}}_{\Bbb {H}}$ be a closed convex cone  generated by   ${\Bbb H}$.
Suppose the random variables  ${\bf X}$  and ${\bf Y}$  are as in (2.5) with finite second moments.
If ${\bf X}\le _ {\cal {F}_{\Bbb{H}}\cap{\Bbb I}} {\bf Y}$, then,
 $\boldsymbol{\mu}^x+ E(|U_0|)\boldsymbol{\delta}_w^x\le \boldsymbol{\mu}^y+ E(|U_0|)\boldsymbol{\delta}_w^y$ and
${\bf \Omega}^{y}-{\bf  \Omega }^x\in {\cal{C}}_{\Bbb H}^*$.
\end{theorem}

  A  random variable ${X}$ is said to have a  skew-elliptical  distribution  with location parameter $\mu$, scale parameter $\sigma$  and skew parameter $\delta$, denoted by
 $X\sim SE_1(\mu, \sigma^2, \delta, g)$,   if its density has the following form
$$f_X(x)=\frac{2}{\sigma}g\left(\frac{(x-\mu)^2}{\sigma^2}\right)G\left( \alpha\frac{x-\mu}{\sigma}\right),$$
where $g$ and $G$  are the density generator and c.d.f. of an elliptical distribution, respectively

The following technical assumption is  needed.

{\bf Assumption 3.1.} We assume density generator $g$ satisfies for  $\sigma_1\neq\sigma_2$,
$$\lim_{t\to-\infty} \frac{\sigma_1}{\sigma_2}\frac{g\left(\frac{t^2}{\sigma^2_2} \right)}{g\left(\frac{t^2}{\sigma^2_1} \right)}
=\lim_{t\to+\infty} \frac{\sigma_1}{\sigma_2}\frac{g\left(\frac{t^2}{\sigma^2_2} \right)}{g\left(\frac{t^2}{\sigma^2_1} \right)}=C,$$
where $C\in[-\infty,1)\cup (1,\infty]$.

Next lemma gives necessary conditions for the usual stochastic ordering of univariate continuous random variables.

\begin{lemma} (Jamali et al. (2020)).  Let $X_1$ and $X_2$ be two random variables with pdfs $f_1$ and $f_1$, respectively. If $X_1\le_{st} X_2$, then
$$\lim_{t\to-\infty} \frac{f_2(t)}{f_1(t)}\le 1\;\; and \;\;  \lim_{t\to+\infty} \frac{f_2(t)}{f_1(t)}\ge 1.$$
\end{lemma}

The following lemma provides the necessary/sufficient conditions for the usual
stochastic ordering of the univariate skew-elliptical distributions. The  univariate skew-normal
case can be found in  Jamali,  Amiri \& Jamalizadeh  (2020).

\begin{lemma}  Let $X\sim SE_1(\mu_x, \sigma_x^2, \delta_x, g)$ and $Y\sim SE_1(\mu_y, \sigma_y^2, \delta_y, g)$ with  $\delta_x \delta_y\ge 0$.

(1) If $\mu_x\le \mu_y$,  $\sigma_x=\sigma_y$ and  $\delta_x\le \delta_y$, then  $X \le_{st} Y$.

(2) If  $X \le_{st} Y$,   then  ${\mu}_x\le {\mu}_y$,  ${\delta}_x\le{\delta}_y$   and $\sigma_x=\sigma_y$,   provided that the density generator $g_1$ satisfies Assumption 3.1.
\end{lemma}
{\bf Proof}  Part (1) follows immediately from Corollary 2.1.

(2) If $X \le_{st} Y$, then $EX\le EY$, and thus   ${\mu}_x+ E(|U_0|)\sigma_x {\delta}_x\le {\mu}_y+ E(|U_0|)\sigma_y{\delta}_y$. By Lemma 3.1, if $X\le_{st} Y$, then
$$\lim_{t\to-\infty} \frac{\sigma_x}{\sigma_y}\frac{g\left(\frac{(t- {\mu}_y)^2}{\sigma^y_2} \right)}{g\left(\frac{(t- {\mu}_x)^2}{\sigma^2_x} \right)}
\frac{G\left(\frac{\alpha_2(t- {\mu}_y)}{\sigma_y} \right)}{G \left(\frac{\alpha_1(t- {\mu}_x)}{\sigma_x} \right)}\le 1,$$
and
$$\lim_{t\to+\infty}   \frac{\sigma_1}{\sigma_2}\frac{g\left(\frac{(t- {\mu}_y)^2}{\sigma^2_2} \right)}{g\left(\frac{(t- {\mu}_x)^2}{\sigma^2_1} \right)}
\frac{G\left(\frac{\alpha_2(t- {\mu}_y)}{\sigma_2} \right)}{G \left(\frac{\alpha_1(t- {\mu}_x)}{\sigma_1} \right)} \ge 1,$$
where $G$ is the CDF of $g$, and
$$\alpha_1=\frac{\delta_x}{\sigma_x\sqrt{1-(\frac{\delta_x}{\sigma_x})^2}},\;\; \alpha_2=\frac{\delta_y}{\sigma_y\sqrt{1-(\frac{\delta_y}{\sigma_y})^2}}.$$
But, under condition  $\delta_x \delta_y\ge 0$, this is only possible if  $\sigma^2_x=\sigma^2_y$.  Besides,   using the same argument as that of  Amiri et al. (2020), the conditions ${\mu}_x\le {\mu}_y$,  ${\delta}_x\le{\delta}_y$ have to hold.
$\hfill\square$

\begin{lemma}  Let $X\sim SE_1(\mu_x, \sigma_y^2, \delta_x, \phi)$ and $Y\sim SE_1(\mu_y, \sigma_y^2, \delta_y, \phi)$.

(1) If $\mu_x\le\mu_y$,  $\sigma_x\le \sigma_y$ and  $\delta_x\le \delta_y$, then  $X \le_{icx} Y$.

(2) If  $X \le_{icx} Y$ and ${\mu}_x={\mu}_y$, then  $\sigma_x {\delta}_x\le \sigma_y{\delta}_y$  and $\sigma_x\le \sigma_y$,  provided that   $X$ and $Y$ are   supported on ${\Bbb{R}}$.
\end{lemma}
{\bf Proof}  Part (1) follows immediately from Corollary 2.1.
(2) If $X \le_{icx} Y$, then $EX\le EY$, and thus    ${\mu}_x+ E(|U_0|)\sigma_x {\delta}_x\le {\mu}_y+ E(|U_0|)\sigma_y{\delta}_y$ which implies
  $\sigma_x {\delta}_x\le \sigma_y{\delta}_y$  since ${\mu}_x={\mu}_y$. Because  $X$ and $Y$ have the same generator, so they have the following stochastic representations:
$$X=\mu_x+ \sigma_x  {\delta}_x|U_0|+  \sigma_x \sqrt{1-{\delta}_x^2}U_1\;\; and \;\; Y=\mu_y+ \sigma_y  {\delta}_y|U_0|+  \sigma_y \sqrt{1-{\delta}_y^2}U_1,$$
where $(U_0,U_1)$.  If  $\sigma_x \sqrt{1-{\delta}_x^2}>\sigma_y \sqrt{1-{\delta}_y^2}$, then  for large $t$,
\begin{eqnarray*}
E(X-t)_+&=&\int_t^{\infty}P\left(U_1>\frac{s-\mu_x- \sigma_x  {\delta}_x|U_0|}{ \sigma_x \sqrt{1-{\delta}_x^2}}\right) ds\\
&&>\int_t^{\infty} P\left(U_1>\frac{s-\mu_y- \sigma_y  {\delta}_y|U_0|}{ \sigma_y \sqrt{1-{\delta}_y^2}}\right)ds=E(Y-t)_+,
\end{eqnarray*}
which contradicts to $X \le_{icx} Y$. Hence, we conclude that
$$\sigma_x \sqrt{1-{\delta}_x^2}\le \sigma_y \sqrt{1-{\delta}_y^2},$$
which implies
 $\sigma_y^2-\sigma_x^2\ge \sigma_y^2\delta_y^2-\sigma_x^2\delta_x^2\ge 0$.
 $\hfill\square$

\begin{lemma}  Let $X\sim SE_1(\mu_x, \sigma_y^2, \delta_x, \phi)$ and $Y\sim SE_1(\mu_y, \sigma_y^2, \delta_y, \phi)$.

(1) If $\mu_x=\mu_y$,  $\sigma_x\le \sigma_y$ and  $\delta_x=\delta_y$, then  $X \le_{cx} Y$.

(2) If  $X \le_{cx} Y$ and ${\mu}_x={\mu}_y$, then  $\sigma_x {\delta}_x=\sigma_y{\delta}_y$  and $\sigma_x\le \sigma_y$,  provided that   $X$ and $Y$ are   supported on ${\Bbb{R}}$.
\end{lemma}
{\bf Proof}  Part (1) follows immediately from Corollary 2.1.  (2) It is well-known that   $X \le_{cx} Y$  if, and only if,  $X \le_{icx} Y$ and $E(X)=E(Y)$. Thus, by Lemma 3.3,   $\sigma_x\le \sigma_y$ and  ${\mu}_x+ E(|U_0|)\sigma_x {\delta}_x={\mu}_y+ E(|U_0|)\sigma_y{\delta}_y$, which implies
  $\sigma_x {\delta}_x=\sigma_y{\delta}_y$  since ${\mu}_x={\mu}_y$.

The following result  gives necessary conditions for   some special cases of increasing-Hessian orders.
\begin{theorem}
Suppose the random variables  ${\bf X}$  and ${\bf Y}$  are as in (2.5).\\
(1) If ${\bf X}\le_{st} {\bf Y}$ with  ${\delta}_{iw}^x{\delta}_{iw}^y\ge 0, i=1,2,\ldots,n$,   then, $\boldsymbol{\mu}^x\le \boldsymbol{\mu}^y$,
 $\boldsymbol{\delta}_w^x\le \boldsymbol{\delta}_w^y$ and ${\bf \Omega}^{y}={\bf  \Omega }^x$;\\
(2) If ${\bf X}\le _{icx} {\bf Y}$ and  $\boldsymbol{\mu}^x=\boldsymbol{\mu}^y$, then,  $\boldsymbol{\delta}_w^x\le \boldsymbol{\delta}_w^y$ and
${\bf\Omega}^{y}-{\bf  \Omega}^x$ is copositive;\\
(3) If  $({\bf X}-E({\bf X}))\le _{ism} ({\bf Y}-E({\bf Y}))$, then,   $\omega_{ii}^x=   \omega_{ii}^y$   for all $1\le i \le n$ and  $\omega_{ij}^x\le \omega_{ij}^y$   for all $1\le i<i \le n$;\\
(4)  If  $({\bf X}-E({\bf X}))\le _{idcx} ({\bf Y}-E({\bf Y}))$, then,   $\omega_{ij}^x\le   \omega_{ij}^y$ for all $1\le i,j\le n$;\\
(5) If ${\bf X}\le _{iccx} {\bf Y}$   and  $\boldsymbol{\mu}^x=\boldsymbol{\mu}^y$, then,  $\boldsymbol{\delta}_w^x\le \boldsymbol{\delta}_w^y$ and  $\omega_{ij}^x\le   \omega_{ij}^y$ for all $1\le i,j\le n$;\\
(6)  If  ${\bf X}\le_{icp} {\bf Y}$, then, $\boldsymbol{\mu}^x+ E(|U_0|)\boldsymbol{\delta}_w^x\le \boldsymbol{\mu}^y+ E(|U_0|)\boldsymbol{\delta}_w^y$ and ${\bf  \Omega}^y-{\bf  \Omega}^x$ is copositive;\\
(7)  If  ${\bf X}\le_{icop} {\bf Y}$,  then, $\boldsymbol{\mu}^x+ E(|U_0|)\boldsymbol{\delta}_w^x\le \boldsymbol{\mu}^y+ E(|U_0|)\boldsymbol{\delta}_w^y$   and ${\bf \Omega}^y-{\bf  \Omega}^x$ is  completely copositive;\\
(8) If ${\bf X}\le _{uo} {\bf Y}$ and  $\boldsymbol{\delta}_w^x=\boldsymbol{\delta}_w^y$,  then, $\boldsymbol{\mu}^x \le \boldsymbol{\mu}^y$,   $\omega_{ii}^x=  \omega_{ii}^y$  for all $1\le i \le n$.
\end{theorem}
{\bf Proof}. 
(1)   We  choose  $f$ to have the forms of
$f({\bf x})=h_1(x_i)$ and $f({\bf x})=h_2(x_i+x_j)$, where  $h_1$ and $h_2$  are any two univariate increasing functions,  it follows from  ${\bf X}\le_{st} {\bf Y}$ that
$X_i\le_{st} Y_i$ and  $X_i+X_j\le_{st} Y_i+Y_j$. Note that  ${\bf X}\sim SE_{n}\left(\boldsymbol{\mu}^x, \mathbf{ \Omega}^x,\boldsymbol{\alpha}_w^x, g^{(n+1)}\right)$ and  ${\bf{Y}}\sim SE_n ({\boldsymbol \mu}^y, \mathbf{ \Omega}^y,\boldsymbol{\alpha}_w^y,  g^{(n+1)})$
 lead to  $X_i\sim SE_1 ({\mu}_i^x,{\omega}_{ii}^x,  {\alpha}_{iw}^x,  g^{(1)})$, $Y_i\sim SE_1 ({\mu}_i^y,{\omega}_{ii}^y,  {\alpha}_{iw}^y,  g^{(1)})$,   $X_i+X_j\sim SE_1 ({\mu}_i^x+{\mu}_j^x,{\omega}_{ii}^x+ {\omega}_{jj}^x+2\omega^x_{ij}, {\alpha}_{iw}^x+ {\alpha}_{jw}^x,  g^{(1)})$ and $Y_i+Y_j\sim SE_1 ({\mu}_i^y+{\mu}_j^y,{\omega}_{ii}^y+ {\omega}_{jj}^y+2\omega^y_{ij},  {\alpha}_{iw}^y+ {\alpha}_{jw}^y,  g^{(1)})$.   Applying Lemma 3.2 (2) we find that  ${\mu_i}^x\le {\mu_i}^y$, ${\delta}_{iw}^x\le {\delta}_{iw}^y$ and
 $\omega^x_{ij}=\omega^y_{ij}$ for all $1\le i,j\le n$. Hence, $\boldsymbol{\mu}^x\le \boldsymbol{\mu}^y$, $\boldsymbol{\delta}_w^x\le \boldsymbol{\delta}_w^y$ and  ${\bf \Omega}^y={\bf\Omega}^x$.

(2) ${\bf X}\le _{icx} {\bf Y}$   implies $\boldsymbol{\mu}^x+ E(|U_0|)\boldsymbol{\delta}_w^x\le \boldsymbol{\mu}^y+ E(|U_0|)\boldsymbol{\delta}_w^y$, this together with  $\boldsymbol{\mu}^x=\boldsymbol{\mu}^y$ yield  $\boldsymbol{\delta}_w^x\le \boldsymbol{\delta}_w^y$.
  For all ${\bf a}\ge 0$,  the function $f_a({\bf x})=f({\bf a}'{\bf x})$ is increasing convex for all increasing convex functions $f: \Bbb{R} \rightarrow \Bbb{R}$.
Hence ${\bf X}\le _{icx} {\bf Y}$   implies
  ${\bf a}'{\bf X}\le _{icx} {\bf a}'{\bf Y}$ for ${\bf a\ge 0}$. Note that ${\bf a}'{\bf X}\sim SE_1({\bf a}'\boldsymbol{\mu}^x, {\bf a}'{\Omega}^x{\bf a}, {\bf a}'\boldsymbol{\delta}_w^x, g^{(1)})$ and
  ${\bf a}'{\bf Y}\sim SE_1({\bf a}'\boldsymbol{\mu}^y, {\bf a}'{\Omega}^y{\bf a}, {\bf a}'\boldsymbol{\delta}_w^y, g^{(1)})$,
  it then follows from Lemma 3.3 (2) that
 ${\bf a}'({\bf\Omega}^{y}-{\bf  \Omega}^x){\bf a}={\bf a}'{\bf\Omega}^{y}{\bf a}-{\bf a}'{\bf  \Omega}^x{\bf a} \ge 0$.

(3)   $({\bf X}-E({\bf X}))\le _{ism} ({\bf Y}-E({\bf Y}))$  implies  $({X}_i-E(X_i))\le _{st} ({Y}_i-E(Y_i))$, and hence $({X}_i-E(X_i))\overset{d}{=}({Y}_i-E(Y_i))$
since they have the same mean. Thus, $\omega_{ii}^x=   \omega_{ii}^y$   for all $1\le i \le n$. Consequently, $({\bf X}-E({\bf X}))\le _{sm} ({\bf Y}-E({\bf Y}))$. Applying  Theorem 3.4(b), we get
$ \omega_{ij}^x\le  \omega _{ij}^y$ for all $1\le i\neq j\le n$.

(4)   $({\bf X}-E({\bf X}))\le _{idcx} ({\bf Y}-E({\bf Y}))$  implies   $({\bf X}-E({\bf X}))\le _{dcx} ({\bf Y}-E({\bf Y}))$ since they have the same mean. Applying  Theorem 3.4(c), we get $ \omega_{ij}^x\le  \omega _{ij}^y$ for all $1\le i, j\le n$.

(5)  If ${\bf X}\le _{iccx} {\bf Y}$, then,  ${\bf X}\le _{icx} {\bf Y}$ (see M\"uller and Stoyan (2002)).   By (2) we get $\boldsymbol{\delta}_w^x\le \boldsymbol{\delta}_w^y$ and ${\Omega}^y-{\bf  \Omega}^x$ is copositive. The latter inequality implies that  $\omega_{ij}^x\le   \omega_{ij}^y$ for all $1\le i,j\le n$.

(6)  If  ${\bf X}\le_{icp} {\bf Y}$, then, ${\bf a}'{\bf X}\le_{icx} {\bf a}'{\bf Y}$ for $a\in{\Bbb{R}}^n_+$ (see Amiri et al.(2022),p 693), similar to (2), we get   ${\bf  \Omega}^y-{\bf  \Omega}^x$ is copositive.

(7)  For all ${\bf a}\ge 0$, then, the function $f_a({\bf x})=f({\bf a}'{\bf x})$ is increasing convex for all increasing convex functions $f: \Bbb{R} \rightarrow \Bbb{R}$.
  The Hessian matrix is given as ${\bf H}_f({\bf x})=({\bf a}{\bf a}')f''({\bf a}'{\bf x})$. Since  $f$ is convex, we have $f''({\bf a}'{\bf x})\ge 0$.
 Hence  ${\bf H}_f({\bf x})$ is completely positive.   If  ${\bf X}\le_{icop} {\bf Y}$,  then,  ${\bf a}'{\bf X}\le_{icx} {\bf a}'{\bf Y}$ for $a\in{\Bbb{R}}^n_+$, similar to (2), we get   ${\bf  \Omega}^y-{\bf  \Omega}^x$ is copositive.

(8)  ${\bf X}\le _{uo} {\bf Y}$ implies $\boldsymbol{\mu}^x+ E(|U_0|)\boldsymbol{\delta}_w^x\le \boldsymbol{\mu}^y+ E(|U_0|)\boldsymbol{\delta}_w^y$ and  ${X}_i\le _{st} {Y}_i$ for all $1\le i\le  n$, these together with  $\boldsymbol{\delta}_w^x=\boldsymbol{\delta}_w^y$ lead to $\boldsymbol{\mu}^x \le \boldsymbol{\mu}^y$,
$\omega_{ii}^x=\omega_{ii}^y$   for all $1\le i \le n$. $\hfill\square$

\subsection{Linear orderings }\label{intro}

\begin{theorem}
Suppose the random variables  ${\bf X}$  and ${\bf Y}$  are as in (2.5).\\
(1)  If ${\boldsymbol \mu}^x\le {\boldsymbol \mu}^y, {\boldsymbol \delta}_w^x\le {\boldsymbol \delta}_w^y$ and ${\bf  \Omega}^{x} = {\bf \Omega}^y$, then ${\bf X}\le _{plst} {\bf Y}$;\\
(2) If ${\bf X}\le_{plst} {\bf Y}$  and ${\delta}_{iw}^x{\delta}_{iw}^y\ge 0, i=1,2,\ldots,n$,  then, $\boldsymbol{\mu}^x\le \boldsymbol{\mu}^y$,
 $\boldsymbol{\delta}_w^x\le \boldsymbol{\delta}_w^y$ and ${\bf \Omega}^{y}={\bf  \Omega }^x$.
\end{theorem}
{\bf Proof}. (1). By Theorem 3.6(1),  if ${\boldsymbol \mu}^x\le {\boldsymbol \mu}^y, {\boldsymbol \delta}_w^x\le {\boldsymbol \delta}_w^y$ and ${\bf  \Omega}^{x} = {\bf \Omega}^y$, then ${\bf X}\le _{st} {\bf Y}$, so we have ${\bf X}\le _{plst} {\bf Y}$.\\
(2) We note that  ${\bf X}\le_{plst} {\bf Y}$ implies $X_i\le _{st} Y_i$ and   $X_i+X_j\le _{st} Y_i+Y_j$ for all $i,j\in\{1,2,\ldots,n\}$. Using Lemma 3.2, together with ${\delta}_{iw}^x{\delta}_{iw}^y\ge 0$, yield  $\boldsymbol{\mu}^x\le \boldsymbol{\mu}^y$,    $\boldsymbol{\delta}_w^x\le \boldsymbol{\delta}_w^y$ and ${\bf \Omega}^{y}={\bf  \Omega }^x$. $\hfill\square$

\begin{theorem}
Suppose the random variables  ${\bf X}$  and ${\bf Y}$  are as in (2.5).\\
(1)   If ${\boldsymbol \mu}^x= {\boldsymbol \mu}^y,   {\boldsymbol \delta}_w^x={\boldsymbol \delta}_w^y$ and ${\bf\Omega}^{y}-{\bf\Omega}^x\succeq O$, then ${\bf X}\le _{lcx} {\bf Y}$;\\
(2) If ${\bf X}\le_{lcx} {\bf Y}$,    then,  $\boldsymbol{\mu}^x+ E(|U_0|)\boldsymbol{\delta}_w^x=\boldsymbol{\mu}^y+ E(|U_0|)\boldsymbol{\delta}_w^y$
  and  ${\bf\Omega}^{y}-{\bf\Omega}^x\succeq O$.
\end{theorem}
{\bf Proof}. (1).  If ${\boldsymbol \mu}^x= {\boldsymbol \mu}^y,   {\boldsymbol \delta}_w^x={\boldsymbol \delta}_w^y$ and ${\bf\Omega}^{y}-{\bf\Omega}^x\succeq O$, then, by Theorem 3.2(a),  ${\bf X}\le _{cx} {\bf Y}$, and so that ${\bf X}\le _{lcx} {\bf Y}$.

(2) If ${\bf X}\le_{lcx} {\bf Y}$, then ${\bf a}'{\bf X}\le_{cx}{\bf a}'{\bf Y}$ for any ${\bf a}\in\Bbb{R}^n$.
 Note that ${\bf a}'{\bf X}\sim SE_1({\bf a}'\boldsymbol{\mu}^x, {\bf a}'{\Omega}^x{\bf a}, {\bf a}'\boldsymbol{\delta}_w^x, g^{(1)})$ and
  ${\bf a}'{\bf Y}\sim SE_1({\bf a}'\boldsymbol{\mu}^y, {\bf a}'{\Omega}^y{\bf a}, {\bf a}'\boldsymbol{\delta}_w^y, g^{(1)})$,
  it is well-known that  ${\bf a}'{\bf X}\le_{cx}{\bf a}'{\bf Y}$ if, and only if, $E({\bf a}'{\bf X}) = E({\bf a}'{\bf Y})$ and   ${\bf a}'{\bf X}\le_{icx}{\bf a}'{\bf Y}$.
 It then $E({\bf X}) = E({\bf Y})$ and follows from Lemma 3.3 (2) that
 ${\bf a}'({\bf\Omega}^{y}-{\bf  \Omega}^x){\bf a}={\bf a}'{\bf\Omega}^{y}{\bf a}-{\bf a}'{\bf  \Omega}^x{\bf a} \ge 0$. Hence,  $\boldsymbol{\mu}^x+ E(|U_0|)\boldsymbol{\delta}_w^x=\boldsymbol{\mu}^y+ E(|U_0|)\boldsymbol{\delta}_w^y$
  and  ${\bf\Omega}^{y}-{\bf\Omega}^x\succeq O$. $\hfill\square$

\begin{theorem}
Suppose the random variables  ${\bf X}$  and ${\bf Y}$  are as in (2.5).\\
(1) If ${\boldsymbol \mu}^x= {\boldsymbol \mu}^y, {\boldsymbol \delta}_w^x= {\boldsymbol \delta}_w^y$ and ${\bf\Omega}^{y}-{\bf\Omega}^x$ is copositive, then ${\bf X}\le _{plcx} {\bf Y}$;\\
(2) If ${\bf X}\le _{plcx} {\bf Y}$ and ${\boldsymbol \mu}^x= {\boldsymbol \mu}^y$, then ${\boldsymbol \delta}_w^x= {\boldsymbol \delta}_w^y$ and ${\bf\Omega}^{y}-{\bf\Omega}^x$ is copositive,  provided that    ${\bf X}$ and  ${\bf Y}$ are   supported on ${\Bbb{R}}$.
\end{theorem}
{\bf Proof}. (1).  Assume that $f: \Bbb{R}\rightarrow \Bbb{R}$ is a convex function and consider the function $u({\bf x}):=f({\bf a}'{\bf x})$ for any  $a\in{\Bbb{R}}^n_+$.
Suppose that ${\boldsymbol \mu}^x= {\boldsymbol \mu}^y, {\boldsymbol \delta}_w^x= {\boldsymbol \delta}_w^y$ and ${\bf  \Omega}^{y}-{\bf \Omega}^x$ is copositive, then, it follows from Lemma 2.2 that $Ef({\bf a}'{\bf X})\le  f({\bf a}'{\bf Y})$, since tr$(H_u(\bf x)({\bf\Omega}^{y}-{\bf\Omega}^x))=f''({\bf a}'{\bf x}){\bf a}'({\bf\Omega}^{y}-{\bf\Omega}^x) {\bf a}\ge 0$, which proves ${\bf X}\le _{plcx} {\bf Y}$.

(2) Suppose  ${\bf X}\le _{plcx} {\bf Y}$, then, by definition,  ${\bf a}'{\bf X}\le_{cx} {\bf a}'{\bf Y}$ for all  $a\in{\Bbb{R}}^n_+$ yielding that $X_i\le_{cx} Y_i$ for $i=1,2,\ldots, n $. Hence,  by Lemma 3.4,   ${\boldsymbol \delta}_w^x= {\boldsymbol \delta}_w^y$,  $Var({\bf a}'{\bf X})\le Var({\bf a}'{\bf Y})$ for all  $a\in{\Bbb{R}}^n_+$, which implies ${\bf a}'({\bf\Omega}^{y}-{\bf\Omega}^x) {\bf a}\ge 0$. $\hfill\square$

\begin{theorem}
Suppose the random variables  ${\bf X}$  and ${\bf Y}$  are as in (2.5).\\
(1) If ${\boldsymbol \mu}^x\le {\boldsymbol \mu}^y, {\boldsymbol \delta}_w^x\le {\boldsymbol \delta}_w^y$ and ${\bf\Omega}^{y}-{\bf\Omega}^x$ is copositive, then ${\bf X}\le _{iplcx} {\bf Y}$;\\
(2) If ${\bf X}\le _{iplcx} {\bf Y}$ and ${\boldsymbol \mu}^x= {\boldsymbol \mu}^y$, then ${\boldsymbol \delta}_w^x\le {\boldsymbol \delta}_w^y$ and ${\bf\Omega}^{y}-{\bf\Omega}^x$ is copositive, provided that    ${\bf X}$ and  ${\bf Y}$ are   supported on ${\Bbb{R}}$.
\end{theorem}
{\bf Proof}. (1).  Assume that $f: \Bbb{R}\rightarrow \Bbb{R}$ is an increasing convex function and consider the function $u({\bf x}):=f({\bf a}'{\bf x})$ for any  $a\in{\Bbb{R}}^n_+$.
Suppose that ${\boldsymbol \mu}^x\le {\boldsymbol \mu}^y, {\boldsymbol \delta}_w^x\le {\boldsymbol \delta}_w^y$ and ${\bf  \Omega}^{y}-{\bf \Omega}^x$ is copositive. Then, it follows from Lemma 2.2 that $Ef({\bf a}'{\bf X})\le  f({\bf a}'{\bf Y})$, since $\nabla f({\bf x})\ge 0$ and tr$(H_u(\bf x)({\bf\Omega}^{y}-{\bf\Omega}^x))=f''({\bf a}'{\bf x}){\bf a}'({\bf\Omega}^{y}-{\bf\Omega}^x) {\bf a}\ge 0$, which proves ${\bf X}\le _{iplcx} {\bf Y}$.

(2) Suppose  ${\bf X}\le _{iplcx} {\bf Y}$, then, by definition,  ${\bf a}'{\bf X}\le_{icx} {\bf a}'{\bf Y}$ for all  $a\in{\Bbb{R}}^n_+$ yielding that $X_i\le_{icx} Y_i$ for $i=1,2,\ldots, n $.   Note that
$${\bf a}'{\bf X}\sim SE_1({\bf a}'\boldsymbol{\mu}^x, {\bf a}'{\Omega}^x{\bf a}, {\bf a}'\boldsymbol{\delta}_w^x, g^{(1)})$$
 and
  $${\bf a}'{\bf Y}\sim SE_1({\bf a}'\boldsymbol{\mu}^y, {\bf a}'{\Omega}^y{\bf a}, {\bf a}'\boldsymbol{\delta}_w^y, g^{(1)}).$$
    Hence,  by Lemma 3.3,   ${\boldsymbol \delta}_w^x\le{\boldsymbol \delta}_w^y$ and ${\bf a}'({\bf\Omega}^{y}-{\bf\Omega}^x) {\bf a}\ge 0$. $\hfill\square$

\section{Concluding remarks}\label{intro}

In this paper, we have discussed the Hessian and increasing-Hessian orderings in the skew-elliptical family of distributions. The necessary and/or sufficient conditions for comparing vectors of multivariate skew-elliptical family of distributions have been established
for different stochastic orderings   and some of
their well-known special cases.  Furthermore, the linear forms of usual stochastic, convex and increasing convex,   positive convex  and increasing-positive-convex  orderings
are also discussed. \\

\noindent{\bf Acknowledgements.}    
This research   was supported by the National Natural Science Foundation of China (No. 12071251).\\

\bibliographystyle{model1-num-names}

\end{document}